\documentclass[12pt,oneside,reqno]{amsart}

\usepackage{amssymb}
\usepackage{txfonts}
\usepackage{bbm}
\usepackage{cases}
\usepackage{amsmath}
\usepackage{graphicx}
\usepackage{mathrsfs}
\usepackage{stmaryrd}
\usepackage{amsfonts}
\usepackage{enumerate,amsmath,amssymb,amsthm}

\numberwithin{equation}{section}

\newcommand{\be}{\begin{eqnarray}}
\newcommand{\ee}{\end{eqnarray}}
\newcommand{\ce}{\begin{eqnarray*}}
\newcommand{\de}{\end{eqnarray*}}
\newtheorem{theorem}{Theorem}[section]
\newtheorem{lemma}[theorem]{Lemma}
\newtheorem{remark}[theorem]{Remark}
\newtheorem{definition}[theorem]{Definition}
\newtheorem{proposition}[theorem]{Proposition}
\newtheorem{Examples}[theorem]{Example}
\newtheorem{corollary}[theorem]{Corollary}

\def\eps{\varepsilon}

\def\e{\mathrm{e}}

\def\[{{\Big[}}
\def\]{{\Big]}}
\def\<{{\langle}}
\def\>{{\rangle}}
\def\({{\Big(}}
\def\){{\Big)}}

\def\bx{{\mathbf{x}}}

\def\dif{{\mathord{{\rm d}}}}

\def\={&\!\!=\!\!&}

\def\cJ{{\mathcal J}}

\def\mB{{\mathbb B}}

\def\mE{{\mathbb E}}

\def\mH{{\mathbb H}}

\def\mL{{\mathbb L}}

\def\mN{{\mathbb N}}

\def\mP{{\mathbb P}}

\def\mR{{\mathbb R}}

\def\mX{{\mathbb X}}

\def\1{{\mathbf{1}}}

\def\sA{{\mathscr A}}
\def\sB{{\mathscr B}}

\def\sF{{\mathscr F}}
\def\sG{{\mathscr G}}

\def\sI{{\mathscr I}}

\def\sM{{\mathscr M}}
\def\sN{{\mathscr N}}

\def\sR{{\mathscr R}}
\def\sS{{\mathscr S}}

\def\sX{{\mathscr X}}

\def\geq{\geqslant}
\def\leq{\leqslant}

\def\eps{\varepsilon}

\def\e{\mathrm{e}}

\def\[{{\Big[}}
\def\]{{\Big]}}
\def\<{{\langle}}
\def\>{{\rangle}}
\def\({{\Big(}}
\def\){{\Big)}}

\def\bx{{\mathbf{x}}}

\def\dif{{\mathord{{\rm d}}}}

\def\={&\!\!=\!\!&}
\def\bt{\begin{theorem}}
\def\et{\end{theorem}}
\def\bl{\begin{lemma}}
\def\el{\end{lemma}}
\def\br{\begin{remark}}
\def\er{\end{remark}}
\def\bx{\begin{Examples}}
\def\ex{\end{Examples}}
\def\bd{\begin{definition}}
\def\ed{\end{definition}}
\def\bp{\begin{proposition}}
\def\ep{\end{proposition}}
\def\bc{\begin{corollary}}
\def\ec{\end{corollary}}

\def\geq{\geqslant}
\def\leq{\leqslant}

\def\<{\langle} \def\>{\rangle}

 \def\beq{\begin{equation}}  
 
\def\e{\text{\rm{e}}}

\allowdisplaybreaks

\begin{document}

\title{Stochastic integrals and BDG's inequalities in Orlicz-type spaces}

\author{YINGCHAO XIE and XICHENG ZHANG}

\date{}

\thanks{
The research of YX is partially supported by NSFs of China (11271169) 
and Project Funded by the Priority Academic Program Development (PAPD) of Jiangsu Higher Education Institutions.
The research of XZ is partially supported by NNSFC grant of China (Nos. 11271294, 11325105).}

\address{
Yingchao Xie: School of Mathematics and Statistics,
Jiangsu Normal University Xuzhou, Jiangsu, 221116, P.R. China\\
Email: ycxie@jsnu.edu.cn
 }

\address{
Xicheng Zhang: School of Mathematics and Statistics,
Wuhan University, Wuhan, Hubei 430072, P.R.China\\
Email: XichengZhang@gmail.com
 }

\begin{abstract}
In this paper we extend an inequality of Lenglart, L\'epingle and Pratelli \cite[Lemma 1.1]{LLP} to
general continuous adapted stochastic processes with values in topology spaces.
By this inequality we show Burkholder-Davies-Gundy's inequality 
for stochastic integrals  in Orlicz-type spaces (a class of quasi-Banach spaces) with respect to cylindrical Brownian motions.
\end{abstract}

\maketitle \rm

\section{Introduction}

Let $(\Omega,\sF,(\sF_t)_{t\geq 0},\mP)$ be a stochastic basis satisfying the usual conditions, which will be fixed below.
Let $M$ be a continuous local martingale with starting point zero, and $\<M\>$ the quadratic variation process of $M$. 
The celebrated Burkholder-Davies-Gundy inequality (abbreviated as BDG's inequality) 
states that for any $p\in(0,\infty)$ and stopping time $\tau$,
\begin{align}\label{BDG}
\mE\left(\sup_{t\in[0,\tau]}|M_t|^p\right)\asymp_p\mE\left(\<M\>^{p/2}_\tau\right),
\end{align}
where $\asymp_p$ means that both sides are comparable up to a positive constant depending only on $p$.
This inequality is a basic tool in stochastic analysis (see \cite{Ik-Wa} and \cite{Re-Yo}), 
and has a deep connection with harmonic analysis (see \cite{Du}, \cite{Ba-Da}  and references therein).

\vspace{2mm}

We would like to emphasize that the key point in the original Burkholder and Gundy's paper \cite[Theorem 3.1]{Bu-Gu}
is to show some ``good $\lambda$-inequalities" for the martingale transform of martingale difference sequences. According to
\cite[Definition 4.8, p.164]{Re-Yo}, an ordered pair $(X,Y)$ of positive random variables is said to satisfy the ``good $\lambda$"-inequality if
for some $\beta>1$ and each $\lambda>0,\delta\in(0,1)$, there is a constant $c_\delta$ with $c_\delta\to 0$ as $\delta\to 0$ such that
\begin{align}\label{UU1}
\mP(X\geq\beta\lambda; Y<\delta \lambda)\leq c_\delta \mP(X\geq\lambda).
\end{align}
This type of ``good $\lambda$''-inequality immediately implies that for any $p\in(0,\infty)$ and some constant $C=C(p,\beta, c_\delta)>0$,
\begin{align}\label{AA08}
\mE X^p\leq C\mE Y^p.
\end{align}
In fact, by \eqref{UU1}, we have
\begin{align*}
\mP(X\geq\beta\lambda)
&\leq c_\delta\mP(X\geq\lambda)+\mP(Y\geq\delta\lambda).
\end{align*}
Integration both sides with respect to $p\lambda^{p-1}\dif \lambda$ on $[0,\infty)$ yields
$$
\beta^{-p}\mE X^p\leq c_\delta\mE X^p+\delta^{-p}\mE Y^p.
$$
Choosing $\delta$ small enough, we obtain \eqref{AA08}.

\vspace{2mm}

One of the proofs of \eqref{BDG} is based on the following ``good $\lambda$''-inequalities for Brownian motion $B$: 
For any stopping time $\tau$ and $\beta>1, \delta, \lambda>0$, it holds that  (see \cite[p152-154]{Du})
\begin{align}\label{Good}
\begin{split}
\mP(B^*_\tau>\beta\lambda, \tau^{1/2}<\delta\lambda)&\leq\tfrac{\delta^2}{(\beta-1)^2}\mP(B^*_\tau>\lambda),\\
\mP(\tau^{1/2}>\beta\lambda, B^*_\tau<\delta\lambda)&\leq\tfrac{\delta^2}{\beta^2-1}\mP(\tau^{1/2}>\lambda),
\end{split}
\end{align}
where, for a function $f:[0,\infty)\to\mR$, we have used the notation 
$$
f^*_t:=\sup_{s\in[0,t]}|f_s|.
$$
From these estimates and by \eqref{AA08}, we get
$$
\mE|B^*_\tau|^p\asymp_p\mE\tau^{p/2}.
$$
 For general continuous local martingale, BDG's inequality \eqref{BDG} follows by the time-change argument 
(for example, see \cite[Theorem 1.10, p183]{Re-Yo}).
Notice that there is another way of proving \eqref{BDG} by using It\^o's formula (see \cite[Theorem 3.1, p110]{Ik-Wa}). Also, for $p>1$, 
BDG's inequality \eqref{BDG}  can be derived from general discrete martingale version of BDG's inequality 
through discretizing a continuous martingale (see \cite[Theorem 7.2.6]{St}).

\vspace{2mm}

Up to now, there are numerous works to study BDG's type inequalities 
(see \cite{Le}, \cite{LLP}, \cite{BY}, \cite{NVW}, \cite{Co-Ve}, \cite{Ma-Ro} and references therein). 
For the ``good $\lambda$-inequality'' method, we would like to quote one sentence from 
\cite[p.166]{Re-Yo} that ``This method is actually extremely powerful and,
to our knowledge, can be used to prove all the BDG-type inequalities for continuous processes.'' In particular, the following 
Lenglart-L\'epingle-Pratelli's inequality is a quite useful tool  (see \cite[Lemma 1.1]{LLP} and \cite[Lemma 4.1]{BY}):
Let $M$ and $N$ be two {\it increasing} predictable processes with $M_0=N_0=0$. Suppose that there exist $q,\kappa>0$ such that for all stopping times
$\tau'\leq\tau$,
\begin{align}\label{Ep}
\mE (M_\tau-M_{\tau'})^q\leq \kappa\|N_\tau\|^q_\infty \mP(\tau'\leq\tau),
\end{align}
where $\|\cdot\|_\infty$ denotes the supremum norm of a random variable. 
Then for any moderately increasing function $\Phi$ (see Definition \ref{Def0} below), there is a constant $C=C(q,\kappa,\Phi)>0$ such that
$$
\mE\Phi(M_\infty)\leq C\mE\Phi(N_\infty).
$$
A said above, the method of the proof in \cite{LLP} is to show the ``good $\lambda$"-inequality for $(M_\infty, N_\infty)$ 
by \eqref{Ep} and suitable choices of stopping times.
Notice that this inequality also plays a crucial role in the study of local times in \cite{BY}.
In this work we shall prove an ``abstract'' version for the above inequality for general continuous adapted stochastic processes with values in topology spaces.
Such an extension will provide an easy tool for proving BDG's inequalities both for finite and infinite dimensional stochastic integrals such as in UMD spaces 
(cf. \cite{NVW}).

\vspace{2mm}

To state our results, we first introduce some notions. 
\bd\label{Def2}
Let $E$ be a topology space and $\gamma\geq 1$. A bi-continuous function $\rho: E\times E\to[0,\infty)$ is called a $\gamma$-quasi-metric on $E$ if
\begin{enumerate}[(i)]
\item $\rho(x,y)=0$  if and only if $x=y$.\quad (ii) $\rho(x,y)=\rho(y,x)$.
\item [(iii)] $\rho(x,y)\leq \gamma(\rho(x,z)+\rho(y,z))$ for all $x,y,z\in E$.
\end{enumerate}
\ed
\bd\label{Def0}
A function $\Phi:(0,\infty)\to(0,\infty)$ is called moderately increasing if $\Phi$ is a nondecreasing function with $\lim_{s\to 0}\Phi(s)=0$ 
and for some $\lambda>1$ (and so, all $\lambda>0$) and $c_\lambda>0$,
$$
\Phi(\lambda t)\leq c_\lambda\Phi(t),\  t\geq 0.
$$
All the  moderately increasing function is denoted by $\sA_0$.
\ed

Clearly, the space $\sA_0$ is closed under usual composition, multiplication and addition. 
For $\Phi\in\sA_0$, let $\Phi^+(t):=\lim_{s\downarrow t}\Phi(s)$ and $\Phi^-(t):=\lim_{s\uparrow t}\Phi(s)$. Observe that $\Phi^+,\Phi^-$ are still in $\sA_0$.
Our first main result is
\bt\label{TH10}
Let $E$ be a topology space endowed with  a $\gamma$-quasi-metric $\rho$.
Let $\xi$ be an $E$-valued continuous adapted process and $N$ a nonnegative continuous adapted increasing process with $N_0=0$. 
Suppose that there exist $q, \kappa>0$ such that for all bounded stopping times $\tau'\leq \tau$,
\begin{align}\label{Con}
\mE\rho(\xi_{\tau},\xi_{\tau'})^q\leq\kappa\,\|N_\tau\|^q_\infty\mP(\tau'<\tau).
\end{align}
Then for any $\Phi\in\sA_0$, there is a positive constant $C=C(\gamma,q,\kappa,\Phi)$ such that for any bounded stopping time $\tau>0$,
\begin{align}\label{M1}
\mE\left(\sup_{t\in[0,\tau]}\Phi\Big(\rho(\xi_t,\xi_0)\Big)\right)\leq C\mE\Phi (N_\tau).
\end{align}
\et
\begin{proof}
Let $M_t:=\rho(\xi_t, \xi_0)$. Clearly, $M$ is a nonnegative continuous adapted processes with $M_0=0$.
For $\alpha,\delta\geq 0$, we define two stopping times
$$
\sigma^M_\alpha:=\inf\Big\{t\geq 0: M_t\geq \alpha\Big\},\  \ \ \sigma^N_\delta:=\inf\Big\{t\geq 0: N_t>\delta\Big\}.
$$
Here we have used the convention $\inf\{\emptyset\}=\infty$.
Below, we fix a bounded stopping time $\tau>0$, $0\leq\beta<\alpha$ and $\delta>0$. 
Observe that 
$$
M^*_\tau\geq\alpha\mbox{ and }N_\tau\leq\delta\Rightarrow \sigma^M_\alpha\leq \tau\leq\sigma^N_\delta,
$$
and since $\rho$ is a $\gamma$-quasi-metric and $\sigma^M_\beta\leq\sigma^M_\alpha$, which further implies that
\begin{align*}
&\rho\big(\xi_{\tau\wedge\sigma^N_\delta\wedge\sigma^M_\alpha},\xi_{\tau\wedge\sigma^N_\delta\wedge\sigma^M_\beta}\big)
=\rho\big(\xi_{\sigma^M_\alpha},\xi_{\sigma^M_\beta}\big)
\geq\tfrac{1}{\gamma}M_{\sigma^M_\alpha}-M_{\sigma^M_\beta}=\tfrac{\alpha}{\gamma}-\beta.
\end{align*}
This means  that
\begin{align}\label{AA1}
\mP\big(M^*_\tau\geq\alpha, N_\tau\leq\delta\big)
\leq\mP\left(\rho\big(\xi_{\tau\wedge\sigma^N_\delta\wedge\sigma^M_\alpha},\xi_{\tau\wedge\sigma^N_\delta\wedge\sigma^M_\beta}\big)
\geq\tfrac{\alpha}{\gamma}-\beta\right).
\end{align}
On the other hand, noticing that 
$$
\tau\wedge\sigma^M_\beta<\tau\wedge\sigma^M_\alpha\Rightarrow\sigma^M_\beta<\tau\Rightarrow M^*_\tau\geq\beta,
$$
by \eqref{Con}, we have
\begin{align}\label{EE1}
\mE\left[\rho\big(\xi_{\tau\wedge\sigma^N_\delta\wedge\sigma^M_\alpha},\xi_{\tau\wedge\sigma^N_\delta\wedge\sigma^M_\beta}\big)^q\right]\leq
 \kappa\delta^q\mP(M^*_\tau\geq\beta).
\end{align}
Therefore, by \eqref{AA1} and Chebyschev's inequality, we obtain
\begin{align*}
\mP\big(M^*_\tau\geq\alpha\big)&\leq\mP\big(M^*_\tau\geq\alpha, N_\tau\leq\delta\big)+\mP\big(N_\tau>\delta\big)\\
&\leq\tfrac{\kappa\delta^q}{(\alpha/\gamma-\beta)^q}\mP(M^*_\tau\geq \beta)+\mP\big(N_\tau>\delta\big).
\end{align*}
In particular, taking $\alpha=\lambda,\beta=\lambda/(2\gamma),\delta=\eps\lambda$, we get for all $\lambda,\eps>0$,
$$
\mP\big(M^*_\tau\geq\lambda\big)\leq\kappa(2\gamma\eps)^q\mP(2\gamma M^*_\tau\geq \lambda)+\mP\big(N_\tau>\eps\lambda\big).
$$
For $\Phi\in\sA_0$, let $\phi(s):=\sup_{t>0}\Phi(st)/\Phi(t)$ and $\Phi_n:=\Phi\wedge n$. 
Integrating both sides with respect to $\dif\Phi^+_n(\lambda)$ on $[0,\infty)$,  we obtain
\begin{align*}
\mE \Phi^+_n(M^*_\tau)&\leq\kappa(2\gamma\eps)^q\mE\Phi^+_n (2\gamma M^*_\tau)+\mE \Phi^-_n (N_\tau/\eps)\\
&\leq \kappa(2\gamma\eps)^q\phi(2\gamma)\mE\Phi^+_n (M^*_\tau)+\phi(1/\eps)\mE \Phi^-_n (N_\tau).
\end{align*}
Choosing $\eps$ small enough and letting $n\to\infty$, we get 
$$
\mE\left(\sup_{t\in[0,\tau]}\Phi^+\Big(\rho(\xi_t,\xi_0)\Big)\right)\leq \phi(1/\eps)\big(1-\kappa(2\gamma\eps)^q\phi(2\gamma)\big)^{-1}\mE\Phi^- (N_\tau).
$$
The desired estimate now follows by the increasing property of $\Phi$.
\end{proof}
The following corollary is quite useful for proving BDG's inequalities for continuous local martingales both in finite and infinite dimensional spaces.
\bc\label{TH1}
Let $E_i, i=1,2$ be two topology spaces endowed with a $\gamma_i$-quasi-metric $\rho_i$ respectively.
Let $\xi$ and $\eta$ be $E_1$ and $E_2$-valued continuous adapted processes, respectively. 
Suppose that there is a constant $\kappa>0$ such that for all bounded stopping times $\tau'\leq \tau$,
\begin{align}\label{Con0}
\mE\rho_1(\xi_{\tau},\xi_{\tau'})\leq\kappa\,\mE\rho_2(\eta_{\tau},\eta_{\tau'}).
\end{align}
Then for any $\Phi\in\sA_0$, there is a positive constant $C=C(\gamma_1,\gamma_2,\kappa,\Phi)$ such that for any bounded stopping time $\tau>0$,
\begin{align}\label{M11}
\mE\left(\sup_{t\in[0,\tau]}\Phi\Big(\rho_1(\xi_t,\xi_0)\Big)\right)\leq C\mE\left(\sup_{t\in[0,\tau]}
\Phi \Big(\rho_2(\eta_t,\eta_0)\Big)\right).
\end{align}
\ec
\begin{proof}
Let $N_t:=\sup_{s\in[0,t]}\rho_2(\eta_s,\eta_0)$. The desired estimate \eqref{M11} follows by
$$
\mE\rho_2(\eta_{\tau},\eta_{\tau'})=\mE\Big(\rho_2(\eta_{\tau},\eta_{\tau'})1_{\{\tau'<\tau\}}\Big)
\leq 2\gamma_2\|N_{\tau}\|_\infty\mP(\tau'<\tau)
$$
and using Theorem \ref{TH10}.
\end{proof}

By this corollary, we can give a direct proof of \eqref{BDG} without using the time-change. In fact, by using stopping time techniques, 
without loss of generality, we may assume that $M$ is a square integrable martingale, that is,
$$
\mE |M_t|^2<\infty,\ \ t>0.
$$
By Doob's optional theorem, for any bounded stopping times $\tau'\leq\tau$,
\begin{align*}
\mE |M_{\tau}-M_{\tau'}|^2=\mE|M_{\tau}|^2-\mE|M_{\tau'}|^2
=\mE\Big(\<M\>_{\tau}-\<M\>_{\tau'}\Big).
\end{align*}
By using \eqref{M11} twice,  we have for any $\Phi\in\sA_0$,
\begin{align}\label{BDG0}
\mE\left(\sup_{t\in[0,\tau]}\Phi\Big(|M_t|^2\Big)\right)\asymp_\Phi
\mE\Phi \big(\<M\>_\tau\big).
\end{align}
Two-sided BDG's inequality \eqref{BDG} now follows by taking $\Phi(t)=t^{p/2}$.

\br
The continuity requirement on $\xi$  in Theorem \ref{TH10} is essential. It seems not liable to extend Theorem \ref{TH10} to c\`adl\`ag processes 
in its current form.
The reason is that BDG's inequality does not hold for discrete martingales with $p\in(0,1)$ (see \cite[Example 8.1]{Bu-Gu}). 
\er

It is interesting that the above elementary inequality allows us to show BDG's inequality for
the stochastic integral  in Orlicz-type spaces (a class of quasi-Banach spaces) with respect to Brownian motions.
It should be emphasized that in \cite{NVW} and \cite{Co-Ve}, the authors have already treated the stochastic integration 
in general UMD spaces and quasi-Banach spaces
in abstract framework under some decoupling assumptions.
If it is not impossible, it seems hard to show that the Orlicz-type space $\mL^\Lambda(\mX;\mH)$ defined below is a UMD space or 
possesses the decoupling properties described in \cite{Co-Ve}. 
Even so, our proof of BDG's inequality is more direct and does not use any abstract functional analysis language.

\

This paper is organized as follows: In the next section, we shall recall some notions about quasi-Banach spaces and Orlicz spaces, and also prove some
basic results for later use. In Section 3, we define the stochastic integral in Orlicz-type spaces and prove BDG's inequalities for the corresponding integrals by Corollary \ref{TH1}.

\section{Quasi-Banach spaces and Orlicz-type spaces}

Let us first recall some notions about quasi-Banach spaces.  Let $\mB$ be a vector space. A quasi-norm $\|\cdot\|_\mB$
is a map from $\mB$ to $[0,\infty)$ with the properties:
\begin{itemize}
\item $\|f\|_\mB=0$ if and only if $f=0$.
\item $\|c f\|_\mB=|c|\cdot\|f\|_\mB$ for all $f\in\mB$ and $c\in\mR$.
\item There is a constant $\gamma_\mB\geq 1$ such that for all $f,g\in\mB$, 
$$
\|f+g\|_\mB\leq \gamma_\mB(\|f\|_\mB+\|g\|_\mB).
$$
\end{itemize}
By a basic theorem due to Aoki \cite{Ao} and Rolewicz \cite{Ro}, for $p=1/(\log_2 \gamma_B +1)$, there is an equivalent $p$-subadditive quasi-norm $\||\cdot\||_\mB$
(cf. \cite{Ka}), that is,
$$
|\|f+g|\|^p_\mB\leq |\|f|\|_\mB^p+|\|g|\|^p_\mB,\ \    f,g\in\mB.
$$
Hence, if we let $d(f,g):=\||f-g\||^p_\mB$, then $d$ is a metric on $\mB$.
We call $(\mB,\|\cdot\|_\mB)$ a quasi-Banach space if $\mB$ is complete with respect to this metric.
More contents about the functional analysis aspect of quasi-Banach spaces are referred to \cite{Ka}.

In this work we shall consider a class of special quasi-Banach spaces. First of all, we introduce a subclass of moderately increasing function space $\sA_0$ as follows:
$$
\sA_1:=\Bigg\{\Lambda\in \sA_0 \mbox{ is continuous and } \lim_{t\to\infty}\Lambda(t)=\infty, \lim_{s\downarrow 0}\sup_{t>0}\Big(\Lambda(st)/\Lambda(t)\Big)=0\Bigg\}.
$$
For example,  define for $\alpha\geq 0$,
$$
\Lambda^\alpha(t):=t^\alpha\big((\log t^{-1})^{-1}\wedge 1\big),\ t>0.
$$
It is easy to see that $\Lambda^\alpha\in \sA_0$ for all $\alpha\geq 0$, and $\Lambda^\alpha\in \sA_1$ only for all $\alpha>0$.

Let $(\mX,\sX,\mu)$ be a $\sigma$-finite measure space, and $(\mB,\|\cdot\|_\mB)$ a quasi-Banach space.
For $\Lambda\in \sA_1$, we introduce a space $\mL^\Lambda(\mX;\mB)$, which consists of
all measurable functions $f: \mX\to \mB$ with 
$$
[f]_{\Lambda}:=\int_\mX \Lambda(\|f(x)\|_\mB)\mu(\dif x)<+\infty,
$$
together with the following Luxemburg-type norm: (cf. \cite{Ad-Fo})
$$
\|f\|_\Lambda:=\inf\big\{\lambda>0: [f/\lambda]_\Lambda\leq 1\big\}.
$$
Observe that for all $f\in\mL^\Lambda(\mX;\mB)$,
\begin{align}\label{AA11}
[f/\|f\|_\Lambda]_\Lambda\leq 1.
\end{align}
In fact, let $\lambda_n$ be a decreasing infimumizing sequence in the definition of $\|f\|_\Lambda$ with limit $\|f\|_\Lambda$. By the monotone convergence theorem, we have
$$
[f/\|f\|_\Lambda]_\Lambda=\lim_{n\to\infty}[f/\lambda_n]_\Lambda\leq 1.
$$
The following proposition shows the relationship between $[f]_\Lambda$ and $\|f\|_\Lambda$.
\bp\label{Pr1} For any $\Lambda\in\sA_1$, there are functions $\phi,\varphi\in\sA_0$ depending only on 
$\Lambda$ with $\lim_{t\uparrow\infty}\phi(t)=\lim_{t\uparrow\infty}\varphi(t)=\infty$
such that for all $f\in\mL^\Lambda(\mX; \mB)$,
\begin{align}\label{AA2}
[f]_\Lambda\leq \phi(\|f\|_\Lambda),\quad \|f\|_\Lambda\leq\varphi([f]_\Lambda).
\end{align}
\ep
\begin{proof}
(1) Let us define
\begin{align}\label{phi}
\phi(s):=\sup_{t>0}\frac{\Lambda(s t)}{\Lambda(t)},\ s\geq 0.
\end{align}
Clearly, $\phi$ is nondecreasing, and $\lim_{s\downarrow 0}\phi(s)=0$ by $\Lambda\in\sA_1$, and
\begin{align}
\phi(st)\leq\phi(s)\phi(t),\ \ s,t\geq 0.\label{EE0}
\end{align}
Hence, $\phi\in\sA_0$ and by definition and \eqref{AA11},
$$
[f]_\Lambda\leq \phi(\|f\|_\Lambda)\int_\mX\Lambda\Big(\|f(x)\|_\mB/\|f\|_\Lambda\Big)\mu(\dif x)\leq\phi(\|f\|_\Lambda).
$$
Moreover, by \eqref{EE0}, we also have
\begin{align*}
1=\phi(1)\leq\phi(1/t)\phi(t)\Rightarrow\lim_{t\to\infty}\phi(t)=\infty.
\end{align*}

(2) To show the existence of $\varphi$, we define 
$$
\psi(t):=\inf\big\{s\geq 0: \phi(s)\geq t\big\},\ \ 
\varphi(t):=1\big/\psi\big(1/t\big).
$$
Since $\phi$ is left continuous and nondecreasing, $\psi$ is also left continuous and nondecreasing, and $\psi(t)>0$ for $t>0$, and
$$
\lim_{t\downarrow 0}\psi(t)=0,\ \lim_{t\to\infty}\psi(t)=\infty,\ \phi(\psi(t))\leq t.
$$
Moreover, by \eqref{EE0} again, for all $\lambda>0$, we also have
\begin{align*}
\psi(t/\phi(\lambda))&=\inf\{s\geq 0: \phi(s)\geq t/\phi(\lambda)\}\\
&\leq\inf\{s\geq 0: \phi(\lambda s)\geq t\}=\psi(t)/\lambda.
\end{align*}
Therefore,  $\varphi\in\sA_0$ and $\lim_{t\to\infty}\varphi(t)=\infty$.
Now by definition, we have
\begin{align*}
[f/\varphi([f]_\Lambda)]_\Lambda&=\int_\mX \Lambda\big(\|f(x)\|_\mB/\varphi([f]_\Lambda)\big)\mu(\dif x)\\
&\leq \phi(1/\varphi([f]_\Lambda))\int_\mX \Lambda(\|f(x)\|_\mB)\mu(\dif x)\\
&=\phi\big(\psi\big(1/[f]_\Lambda\big)\big)[f]_\Lambda\leq1,
\end{align*}
which yields that $\|f\|_\Lambda\leq\varphi([f]_\Lambda)$. The proof is complete.
\end{proof}
\br
From the above constructions of $\phi$ and $\varphi$, it is natural to ask whether $\varphi\circ\phi(t)\leq t$ so that 
$$
\|f\|_\Lambda\leq\varphi([f]_\Lambda)\leq\varphi(\phi(\|f\|_\Lambda))\leq\|f\|_\Lambda.
$$
The answer is in general negative. In fact, assume that $\Lambda\in\sA_1$ is convex. Then, one sees that $\phi(s):=\sup_{t>0}(\Lambda(st)/\Lambda(t))$
is also convex, and hence, continuous and strictly increasing. Thus, $\psi=\phi^{-1}$
and $\varphi(t)=1/\phi^{-1}(1/t)$. If $\varphi\circ\phi(t)\leq t$, then $\phi(t)\phi(1/t)\leq 1$, which together with $\phi(st)\leq\phi(s)\phi(t)$ implies that
$\phi(st)=\phi(s)\phi(t)$. Therefore, $\phi(t)=t^p$ for some $p>0$. This is the case of $L^p$-spaces, but not the case of general Orlicz spaces.
\er

\bp
The space $(\mL^\Lambda(\mX;\mB),\|\cdot\|_\Lambda)$ is a quasi-Banach space, and called Orlicz-type space.
\ep
\begin{proof}
Since $\Lambda$ is moderately increasing, we have for any $c_1,c_2\geq 0$,
\begin{align*}
\Lambda(c_1s+c_2t)&\leq\Lambda(2(c_1\vee c_2)(s\vee t))\leq\phi(2(c_1\vee c_2))\Lambda(s\vee t)\\
&\leq \phi(2(c_1\vee c_2))(\Lambda(s)+\Lambda(t)),
\end{align*}
where $\phi$ is defined by \eqref{phi}.
Hence, 
$$
[c_1 f+c_2 g]_\Lambda\leq\phi(2\gamma_\mB(c_1\vee c_2))([f]_\Lambda+[g]_\Lambda)<\infty,
$$
which implies that $\mL^\Lambda(\mX;\mB)$ is a linear space, and by \eqref{AA11}, for any $\alpha>0$,
$$
\Big[\tfrac{f+g}{\alpha(\|f\|_\Lambda+\|g\|_\Lambda)}\Big]_\Lambda
\leq 2\phi\Big(2\gamma_\mB\tfrac{\|f\|_\Lambda\vee\|g\|_\Lambda}{\alpha(\|f\|_\Lambda+\|g\|_\Lambda)}\Big)
\leq 2\phi(2\gamma_\mB/\alpha).
$$
Letting $\alpha$ be large enough so that $2\phi(2\gamma_\mB/\alpha)\leq 1$ yields that
$$
\|f+g\|_\Lambda\leq \alpha(\|f\|_\Lambda+\|g\|_\Lambda).
$$
Moreover, it is easy to see that $\|c f\|_\Lambda=|c|\cdot\|f\|_\Lambda$ for any $c\in\mR$, and
by \eqref{AA2}, $\|f\|_\Lambda=0$ if and only if $f(x)=0$ for $\mu-$almost all $x\in\mX$. Thus, $\|\cdot\|_\Lambda$ is a quasi-norm. 
The proof of completeness is standard by using \eqref{AA2} (see \cite[Theorem 2.16]{Ad-Fo}). We omit the details.
\end{proof}

Now, we recall the notion of $N$-functions and its complementary function. A function $\Lambda:(0,\infty)\to(0,\infty)$
is called an $N$-function if $\Lambda$ is convex and satisfies $\lim_{t\to 0}\frac{\Lambda(t)}{t}=0$ and
$\lim_{t\to\infty}\frac{\Lambda(t)}{t}=\infty$. All the $N$-functions is denoted by $\sN$. Let $a(t)=\Lambda'_+(t)$ be the right continuous derivative of $\Lambda$.
Then $a$ is nondecreasing and $a(0)=0$, $a(t)>0$ for $t>0$ and $\lim_{t\to\infty}a(t)=\infty$. In particular, $\Lambda(t)=\int^t_0a(s)\dif s$.
Let $\tilde a(t):=\inf\{t: a(s)>t\}$ be the right inverse of $a$. We call $\tilde\Lambda(t):=\int^t_0\tilde a(s)\dif s$ the complementary $N$-function of $\Lambda$.
In particular, the following Young's inequality holds:
$$
st\leq\Lambda(s)+\tilde\Lambda(t),\ s,t\geq 0.
$$

\br
If $\Lambda\in\sA_1\cap\sN$ and $\mB$ is a Banach space, 
then quasi-Banach $(\mL^\Lambda(\mX;\mB),\|\cdot\|_\Lambda)$ becomes a Banach space. When $\mB=\mR$, it is the usual Orlicz space studied in \cite{Kr-Ru}.
\er

To show two-sided BDG's inequality, we introduce the following subclass of $\sA_1$:
\begin{align*}
\sA_2:=\Bigg\{\Lambda&\in\sA_1\cap\sN:  \int^1_0\frac{\Lambda(st)}{s^2}\dif s\leq \kappa_\Lambda\Lambda(t)\mbox{ for some $\kappa_\Lambda>0$.}\Bigg\}
\end{align*}
\bl\label{Le24}
Let $\xi,\eta$ be two nonnegative random variables with
\begin{align}\label{ES1}
\mP(\xi\geq\lambda)\leq\mE(\eta 1_{\{\xi\geq\lambda\}} /\lambda),\ \lambda>0.
\end{align}
For any $\Lambda\in\sA_2$, there is a constant $C=C(\Lambda)>0$ such that
$$
\mE\Lambda(\xi)\leq C\mE\Lambda(\eta).
$$
\el
\begin{proof}
For $\alpha>0$, integrating both sides of \eqref{ES1} with respect to $\dif\Lambda(\alpha\lambda)$ and by the integration by parts formula, we obtain
\begin{align*}
\mE\Lambda(\alpha\xi)&\leq\mE\left(\eta \int^\xi_0 \frac{\dif\Lambda(\alpha\lambda)}{\lambda}\right)
=\mE\left(\eta \left(\frac{\Lambda(\alpha\xi)}{\xi}+\int^{\xi}_0 \frac{\Lambda(\alpha\lambda) }{\lambda^2}\dif\lambda\right)\right)\\
&=\mE\left(\eta \left(\frac{\Lambda(\alpha\xi)}{\xi}+\int^{1}_0 \frac{\Lambda(\alpha\xi\lambda) }{\xi\lambda^2}\dif\lambda\right)\right)
\leq (1+\kappa_\Lambda)\mE\left(\eta \Lambda(\alpha\xi)/\xi\right).
\end{align*}
Noticing that
$$
\tilde\Lambda(\Lambda(t)/t)\leq\Lambda(t),
$$
by Young's inequality, we further have
\begin{align*}
\mE\Lambda(\alpha\xi)&\leq \alpha(1+\kappa_\Lambda)\mE\left(\Lambda(\eta)+\tilde\Lambda( \Lambda(\alpha\xi)/(\alpha\xi))\right)\\
&\leq \alpha(1+\kappa_\Lambda)\left(\mE\Lambda(\eta)+\mE\Lambda(\alpha\xi)\right).
\end{align*}
Letting $\alpha=\frac{1}{2(1+\kappa_\Lambda)}$, we obtain
$$
\mE\Lambda(\xi/(2(1+\kappa_\Lambda))\leq\mE\Lambda(\eta),
$$
which implies the desired estimate by $\Lambda(s)\leq C\Lambda(s/(2(1+\gamma_\Lambda)))$.
\end{proof}

We can extend the classical Doob's maximal inequality as follows:
\bp\label{Pr25}
Let $M_t$ be a c\`adl\`ag martingale. For any $\Lambda\in\sA_2$, there is a constant $C=C(\Lambda)>0$ such that for any bounded stopping time $\tau$,
$$
\mE\left(\sup_{t\in[0,\tau]} \Lambda(|M_t|)\right)\leq C\mE\left(\Lambda(|M_\tau|)\right).
$$
\ep
\begin{proof}
Let $\tau$ be bounded by $T$. By considering stopping martingale $M^\tau_t:=M_{t\wedge\tau}$, without loss of generality, we may assume $\tau=T$.
Let $D$ be a countable dense subset of $[0,T]$ containing the terminal point $T$. Let $D_n\subset D$ be an increasing sequence of finite subsets of $D$ containing $T$
with $\cup_n D_n=D$. 
By \cite[Proposition 1.5, p53]{Re-Yo} and Lemma \ref{Le24}, we have
$$
\mE\left(\sup_{t\in D_n} \Lambda(|M_t|)\right)\leq C\mE\left(\Lambda(|M_T|)\right),
$$
which in turn gives the desired inequality by taking limits $n\to\infty$ and noting 
$\mE\left(\sup_{t\in [0,T]} \Lambda(|M_t|)\right)=\mE\left(\sup_{t\in D} \Lambda(|M_t|)\right)$.
\end{proof}

To have more intuitive pictures about the functions in $\sA_i, i=0,1,2$, we recall the following notion of  regularly varying functions.
\bd\label{Def31}
A measurable function $\phi: (0,\infty)\to (0,\infty)$ is said to vary regularly at zero with index $\alpha\in\mR$ if
$$
\lim_{t\to 0}\frac{\phi(\lambda t)}{\phi(t)}=\lambda^\alpha,\ \ \lambda>0.
$$
We call such $\phi$
a regularly varying function. All regularly varying functions with index $\alpha$ is denoted by $\sR_\alpha$. In particular, the element in $\sR_0$ is called
slowly varying function. Any $\phi\in\sR_\alpha$ can be written as $\phi(t)=t^\alpha\phi_0(t)$ for some $\phi_0\in\sR_0$.
\ed

The following proposition provides useful examples in $\sA_i, i=0,1,2$.
\bp
For $\alpha\geq 0$, let $\sR^+_\alpha$ be the set of all increasing regularly varying functions with $\lim_{t\to 0}\phi(t)=0$
and being bounded away from $0$ and $\infty$ on any compact subset of $(0,\infty)$. Then it holds that
$$
\cup_{\alpha\geq 0}\sR^+_\alpha\subset\sA_0, \ \ \ \cup_{\alpha>0}\sR^+_\alpha\subset\sA_1,\  \ \ \cup_{\alpha>1}\sR^+_\alpha\cap\sN\subset\sA_2.
$$
\ep
\begin{proof}
Let $\alpha\geq 0$ and $\phi\in\sR_\alpha$. By \cite[p.25-28]{Bi-Go-Te}, for any $\delta>0$, there is a constant $C=C(\delta)\geq 1$
such that for all $t,s\in(0,\infty)$,
$$
\frac{\phi(t)}{\phi(s)}\leq C\max\left\{\Big(\frac{t}{s}\Big)^{\alpha+\delta}, \Big(\frac{t}{s}\Big)^{\alpha-\delta}\right\}.
$$
The desired inclusions follow by the above estimate.
\end{proof}


\section{Stochastic integrals in Orlicz-type spaces}

Let $\mH$ be a separable Hilbert space. In this section, we shall define the stochastic integral of 
$\mL^\Lambda(\mX;\mH)$-valued processes with respect to a cylindrical Brownian motion $B$ in $\mH$.
Let $\{B^{(j)}_t, t\geq 0,  j\in\mN\}$ be a sequence of independent one dimensional standard Brownian motion over stochastic basis $(\Omega,\sF,(\sF_t)_{t\geq 0},\mP)$,
and $\{\e_j, j\in\mN\}$ an orthogonormal basis of $\mH$. 
Since every cylindrical Brownian motion $B$ in $\mH$ can be represented as
$$
B_t=\sum_{j\in\mN}B^{(j)}_t\e_j,\ t\geq 0,
$$
without loss of generality, we always assume $\mH=\ell^2$, where $\ell^2$ 
is the usual sequence Hilbert space, that is,
$$
\ell^2=\Bigg\{a=(a^{(1)},a^{(2)}\cdots),\|a\|_{\ell^2}:=\Bigg(\sum_{j\in\mN} |a^{(j)}|^2\Bigg)^{1/2}\Bigg\}.
$$

Below, we introduce the following notations for simplicity.
\begin{itemize}
\item $\Sigma:=\mR_+\times\mX\times\Omega$, $\sG:=\sB(\mR_+)\times\sX\times\sF$.
\item Let $X:\Sigma\to\ell^2$ be a measurable function. Define
$$
\eta^X_t(x,\omega):=\int^t_0 \|X_s(x,\omega)\|^2_{\ell^2}\dif s.
$$
\item For $\Lambda\in\sA_1$ and $t>0$, define
$$
|\|X(\omega)\||_{\Lambda,t}:=\Big[\big(\eta^X_t(\cdot,\omega)\big)^{1/2}\Big]_{\Lambda}.
$$
\end{itemize}
First of all, we introduce the following classes of stochastic processes.
\bd
Let $X:\Sigma\to\ell^2$ be a measurable function and $\Lambda\in\sA_1$.
\begin{enumerate}[(i)]
\item We call $X$ an adapted process if for each $t\geq 0$, $X_t$ is $\sX\times\sF_t/\sB(\ell^2)$-measurable. 
\item  We call $X$ a progressively measurable process if for each $t\geq 0$,
$X\cdot 1_{[0,t]}$ is $\sB([0,t])\times\sX\times\sF_t/\sB(\ell^2)$-measurable. 
\item We call $X$ an elementary process if $X$ takes the following form
\begin{align}\label{For}
X_t(x,\omega)=\sum_{i=0}^{\infty} \xi_i(x,\omega)\cdot 1_{[s_i, s_{i+1})}(t),
\end{align}
where $0=s_0<s_1<\cdots<s_n\uparrow\infty$ and the map $(x,\omega)\mapsto\xi_i(x,\omega)$ is $\sX\times\sF_{s_i}/\sB(\ell^2)$-measurable and satisfies that
$(x,\omega)\mapsto\|\xi_i(x,\omega)\|_{\ell^2}$ is bounded and for all $\omega$, $\xi_i(\cdot,\omega)\in\mL^\Lambda(\mX;\ell^2)$.
All such elementary processes is denoted by $\sS_\Lambda(\Sigma;\ell^2)$.
\item We denote by $\sM_{\Lambda}(\Sigma; \ell^2)$ the space of all progressively measurable processes with
$$
\mP\Big(\omega: |\|X(\omega)\||_{\Lambda,T}<\infty\Big)=1,\ \ \forall T>0.
$$
\end{enumerate}
\ed
\br
By a deep result in stochastic process theory (cf. \cite{De-Me}), any measurable adapted process has a progressively measurable modification.
Thus, in the definition of $\sM_{\Lambda}(\Sigma; \ell^2)$,
one may replace progressively measurable processes with measurable adapted processes.
\er

The following result is purely technical and standard, which states that any $X\in\sM_{\Lambda}(\Sigma; \ell^2)$ can be approximated by elementary processes.
\bp\label{Pr3}
$\sS_\Lambda(\Sigma;\ell^2)$ is dense in $\sM_\Lambda(\Sigma;\ell^2)$ in the sense that for each $X\in\sM_\Lambda(\Sigma;\ell^2)$,
there exists a sequence $X^n\in\sS_\Lambda(\Sigma;\ell^2)$ such that for each $T>0$ and $\omega\in\Omega$,
$$
|\|X^n(\omega)\||_{\Lambda,T}\leq |\|X(\omega)\||_{\Lambda,T},
$$
and
$$
\lim_{n\to\infty}\mE\Big(1\wedge|\|X^n(\cdot)-X(\cdot)\||_{\Lambda,T}\Big)=0.
$$
\ep
\begin{proof}
Let $U_n\subset\mX$ be a sequence of increasing measurable sets with $\mu(U_n)<\infty$ for each $n\in\mN$ and $\cup_n U_n=\mX$. 
Let $\chi_n:[0,\infty)\to[0,1]$ be a continuous function
with $\chi_n(s)=1$ for $s\leq n-1$ and $\chi_n(s)=0$ for $s>n$. Define
$$
X^n_t(x,\omega):=1_{U_n}(x)\ X_t(x,\omega)\ \chi_n(\|X_t(x,\omega)\|_{\ell^2}).
$$
Clearly,  for each $T>0$, we have $|\|X^n(\omega)\||_{\Lambda,T}\leq |\|X(\omega)\||_{\Lambda,T}$, and by the monotone convergence theorem, 
$$
\lim_{n\to\infty}\mE\Big(1\wedge|\|X^n(\cdot)-X(\cdot)\||_{\Lambda,T}\Big)=0.
$$
Next,  for any $n,m\in\mN$, let us define
$$
X^{n,m}_t(x,\omega):=\cJ^m_t(X^n_\cdot(x,\omega)),
$$
where, for a function $f\in L^1_{loc}(\mR_+;\ell^2)$,
$$
\cJ^m_t(f):=\sum_{j=1}^\infty\left(m\int^{j/m}_{(j-1)/m}f_s\dif s\right)1_{[j/m, (j+1)/m)}(t).
$$
From this construction, it is easy to see that $X^{n,m}\in\sS_\Lambda(\Sigma;\ell^2)$ and
$$
\int^T_0 \|X^{n,m}_s(x,\omega)\|^2_{\ell^2}\dif s\leq \int^T_0 \|X^{n}_s(x,\omega)\|^2_{\ell^2}\dif s,\ \ \forall T>0.
$$
By the dominated convergence theorem, for each fixed $n$, we have
$$
\lim_{m\to\infty}\mE\Big(1\wedge |\|X^{n,m}(\cdot)-X^n(\cdot)\||_{\Lambda,T}\Big)=0,\ \ \forall T>0.
$$
Here we have used the following fact: for any $f\in L^2_{loc}(\mR_+;\ell^2)$ and $T>0$,
\begin{align}\label{NH4}
\lim_{m\to\infty}\int^T_0\left\|\cJ^m_s(f)-f_s\right\|^2_{\ell^2}\dif t=0.
\end{align}
Indeed, let $f^n\in C(\mR_+;\ell^2)$ be a sequence of continuous functions with
$$
\lim_{n\to\infty}\int^T_0\left\|f^n_t-f_t\right\|^2_{\ell^2}\dif t=0.
$$
Limit \eqref{NH4} follows by observing 
$$
\int^T_0\left\|\cJ^m_s(f^n-f)\right\|^2_{\ell^2}\dif t\leq\int^T_0\left\|f^n_s-f_s\right\|^2_{\ell^2}\dif t,
$$
and for each fixed $n$,
$$
\lim_{m\to\infty}\int^T_0\left\|\cJ^m_s(f^n)-f^n_s\right\|^2_{\ell^2}\dif t=0.
$$
Finally, the desired sequence of simple processes is obtained by a standard dialganization argument.
\end{proof}

Now we can state  and prove the following main result of this section.
\bt
Let $\Lambda\in\sA_1$. For any $X\in\sM_\Lambda(\Sigma;\ell^2)$, 
there exists a measurable adapted process $\sI^X_t\in\mL^\Lambda(\mX;\mR)$ (called the stochastic integral of $X$ with respect to $B$) 
with the following properties:
\begin{enumerate}[(a)]
\item For $\mu$-almost all $x\in\mX$, it holds that for all $t\geq 0$,
\begin{align}\label{ST}
\int^t_0\|X_s(x)\|^2_{\ell^2}\dif s<\infty\ \ {\rm and}\  \ \sI^X_t(x)=\int^t_0X_s(x)\dif B_s\ a.s.,
\end{align}
where the right hand side is the usual It\^o's integral. In particular, for $\mu$-almost all $x\in\mX$,
the process $t\mapsto \sI^X_t(x,\cdot)$ is a continuous local martingale with square variation process
$\int^t_0\|X_s(x)\|^2_{\ell^2}\dif s$.
\item Let $\Phi\in\sA_0$. For any stopping time $\tau$, we have
\begin{align}\label{AA6}
\mE\left(\Phi\left(\sup_{t\in[0,\tau]}\left[\sI^X_t\right]_{\Lambda}\right)\right)\preceq_{\Lambda,\Phi}
\mE\left(\Phi\left(\Bigg[\Bigg(\int^\tau_0 \|X_s\|^2_{\ell^2}\dif s\Bigg)^{1/2}\Bigg]_{\Lambda}\right)\right),
\end{align}
and if $\Lambda\in\sA_2$, then we also have the reversed inequality
\begin{align}\label{AA61}
\mE\left(\Phi\left(\Bigg[\Bigg(\int^\tau_0 \|X_s\|^2_{\ell^2}\dif s\Bigg)^{1/2}\Bigg]_{\Lambda}\right)\right)\preceq_{\Lambda,\Phi}
\mE\left(\Phi\left(\sup_{t\in[0,\tau]}\left[\sI^X_t\right]_{\Lambda}\right)\right).
\end{align}
Here $A\preceq_{\Lambda,\Phi} B$ means that $A\leq CB$ for some constant $C=C(\Lambda,\Phi)>0$.
\end{enumerate}
\et
\begin{proof}
(i) First of all, let $X\in\sS_\Lambda(\Sigma;\ell^2)$ be an elementary process. The stochastic integral of $X$ with respect to $B$ is naturally defined by
$$
\sI^X_t(x)=\int^t_0 X_s(x)\dif B_s:=\sum_{i=0}^{\infty} \sum_{j\in\mN} \xi^{(j)}_i(x)(B^{(j)}_{t\wedge s_{i+1}}-B^{(j)}_{t\wedge s_i}).
$$
For each $x\in\mX$, noticing that $\xi_i(x,\cdot)\in\sF_{s_i}$ is bounded, one sees that $t\mapsto \sI^X_t(x)$ is a continuous $\mR$-valued square integrable martingale, and
for any bounded stopping time $\tau$,
\begin{align*}
\mE|\sI^X_{\tau}(x)|^2&=\mE\left(\sum_{i=0}^{\infty} \|\xi_i(x)\|^2_{\ell^2}\Big((\tau\wedge s_{i+1})-(\tau\wedge s_i)\Big)\right)\\
&=\mE\left(\int^\tau_0\|X_s(x)\|^2_{\ell^2}\dif s\right)=\mE \eta^X_\tau(x).
\end{align*}
By Doob's optional theorem, for any bounded stopping times $\tau'\leq\tau$,
$$
\mE|\sI^{X}_{\tau}(x)-\sI^{X}_{\tau'}(x)|^2=\mE|\sI^{X}_{\tau}(x)|^2-\mE |\sI^{X}_{\tau'}(x)|^2
=\mE\left(\eta^{X}_{\tau}(x)-\eta^{X}_{\tau'}(x)\right).
$$ 
Fix $\Lambda\in\sA_1$. By Corollary \ref{TH1}, we have for each $x$ and $\tau$,
\begin{align}\label{Eq1}
\mE\Lambda\Bigg(\sup_{s\in[0,\tau]}|\sI^{X}_{s}(x)|\Bigg)\asymp_\Lambda \mE\Lambda\Big(\big(\eta^{X}_\tau(x)\big)^{1/2}\Big).
\end{align}
(ii) Next, let $X\in\sM_\Lambda(\Sigma;\ell^2)$. By Proposition \ref{Pr3}, there is 
a sequence $X^n\in\sS_\Lambda(\Sigma;\ell^2)$ so that for each $T>0$,
\begin{align}\label{NH1}
|\|X^n(\omega)\||_{\Lambda,T}\leq |\|X(\omega)\||_{\Lambda,T},\ \ \forall\omega\in\Omega,
\end{align}
and
\begin{align}\label{NH2}
\lim_{n\to\infty}\mE\Big(1\wedge|\|X^n(\cdot)-X(\cdot)\||_{\Lambda,T}\Big)=0.
\end{align}
For $R>0$, define a stopping time 
$$
\tau_R(\omega):=\inf\Big\{t\geq 0: \||X(\omega)\||_{\Lambda,t}>R\Big\}.
$$
By \eqref{Eq1}, we have for each $x\in\mX$,
$$
\mE\Lambda\Bigg(\sup_{s\in[0,\tau_R]}\big|\sI^{X^n-X^m}_{s}(x)\big|\Bigg)
\asymp_\Lambda \mE\Lambda\Big(\big(\eta^{X^n-X^m}_s(x)\big)^{1/2}\Big).
$$
Integrating both sides with respect to $\mu$ over $\mX$, we get
$$
\mE\Bigg[\sup_{s\in[0,\tau_R]}\big|\sI^{X^n}_{s}-\sI^{X^m}_{s}\big|\Bigg]_\Lambda
\asymp_\Lambda\mE|\|X^n-X^m\||_{\Lambda,\tau_R},
$$
which converges to zero as $n,m\to\infty$ by \eqref{NH1}, \eqref{NH2} and the dominated convergence theorem.
Since $\lim_{R\to\infty}\tau_R=\infty$ almost surely, there exists a measurable adapted process $\sI^X_t(x)$ so that
for $(\mu\times\mP)$-almost all $(x,\omega)$, $t\mapsto \sI^X_t(x,\omega)$ is continuous and
\begin{align}\label{NH5}
\lim_{n\to\infty}\mE\Bigg[\sup_{s\in[0,\tau_R]}\big|\sI^{X^n}_{s}-\sI^{X}_{s}\big|\Bigg]_\Lambda=0,\ \ \forall R>0.
\end{align}
By Fubini's theorem, \eqref{NH2} and \eqref{NH5}, up to extracting a subsequence, 
there is a $\mu$-null set $D_1\subset\mX$ such that for all $x\notin D_1$ and for any $R>0$,
\begin{align}
&\lim_{n\to\infty}\mE\Lambda\left(\Bigg(\int^{\tau_R}_0\|X^n_{s}(x)-X_{s}(x)\|_{\ell^2}^2\dif s\Bigg)^{1/2}\right)=0,\label{NH6}\\
&\quad\lim_{n\to\infty}\mE\Lambda\Bigg(\sup_{s\in[0,\tau_R]}\big|\sI^{X^n}_{s}(x)-\sI^{X}_{s}(x)\big|\Bigg)=0.\label{NH66}
\end{align}
Since $\lim_{t\to\infty}\Lambda(t)=\infty$, there is also a $\mu$-null set $D_2\subset\mX$ such that for all $x\notin D_2$ and any $T>0$,
$$
\mbox{$t\mapsto \sI^X_t(x)$ is continuous adapted and }\int^T_0\|X_s(x)\|^2_{\ell^2}\dif s<\infty,\ a.s.
$$
We now show that for all $x\in (D_1\cup D_2)^c$,  
\begin{align}\label{ER7}
\sI^X_t(x)=\int^t_0X_s(x)\dif B_s,\ \forall t\geq 0.
\end{align}
Fix $x\in (D_1\cup D_2)^c$. For any $R>0$, define a stopping time
$$
\sigma_R(x):=\inf\left\{t>0: \sI^X_t(x)\vee\int^t_0\|X_s(x)\|^2_{\ell^2}\dif s>R\right\}.
$$
Since $\lim_{R\to\infty}\sigma_R(x)=\lim_{R\to\infty}\tau_R=\infty$ almost surely, to show \eqref{ER7}, 
it suffices to prove that
\begin{align}\label{ER8}
\sI^X_{t\wedge\sigma_R(x)\wedge\tau_R}(x)=\int^{t\wedge\sigma_R(x)\wedge\tau_R}_0X_s(x)\dif B_s,\ \forall t\geq 0.
\end{align}
Observe that for each $n\in\mN$,
\begin{align}\label{ER9}
\sI^{X^n}_{t\wedge\sigma_R(x)\wedge\tau_R}(x)=\int^{t\wedge\sigma_R(x)\wedge\tau_R}_0X^n_s(x)\dif B_s,\ \forall t\geq 0.
\end{align}
By BDG's inequality \eqref{BDG0}, we have
\begin{align*}
&\mE\Lambda\left(\sup_{t\in[0,\sigma_R(x)\wedge\tau_R]}\left|\int^t_0 X^n_s(x)\dif B_s-\int^t_0 X_s(x)\dif B_s\right|\right)\\
&\quad\leq \mE\Lambda\left(\left(\int^{\sigma_R(x)\wedge\tau_R}_0 \|X^n_s(x)-X_s(x)\|^2_{\ell^2}\dif s\right)^{1/2}\right),
\end{align*}
which converges to zero as $n\to\infty$ by \eqref{NH6}. By this limit and \eqref{NH66}, \eqref{ER9}, we get \eqref{ER8}.
\\
\\
(iii) Let $\sigma\leq\tau$ be two bounded stopping times and define
$$
X^\sigma_s(x,\omega):=X_s(x,\omega)1_{[\sigma,\infty)}(s).
$$
Clearly, $X^\sigma\in\sM_\Lambda(\Sigma;\ell^2)$. By \eqref{ST} and BDG's inequality \eqref{BDG0} again, we have for $\mu$-almost all $x\in\mX$,
\begin{align}\label{NK1}
\mE\Lambda\Bigg(\sup_{t\in[0,\tau]}\left|\sI^{X^\sigma}_t(x)\right|\Bigg)\asymp_\Lambda 
\mE\Lambda\Big(\big(\eta^{X^\sigma}_\tau(x)\big)^{1/2}\Big).
\end{align}
Observing that
\begin{align}\label{Eq2}
\sI^{X^{\sigma}}_{\tau}(x)=\sI^X_{\tau}(x)-\sI^X_{\sigma}(x),\ \ \eta^{X^{\sigma}}_\tau(x)=\eta^X_{\tau}(x)-\eta^X_{\sigma}(x),
\end{align}
we get
\begin{align}\label{Eq3}
\mE\Lambda\left(|\sI^X_{\tau}(x)-\sI^X_{\sigma}(x)|\right)\preceq_\Lambda\mE\Lambda\left(|\eta^X_{\tau}(x)-\eta^X_{\sigma}(x)|^{1/2}\right).
\end{align}
If $\Lambda\in\sA_2$, by \eqref{NK1}, \eqref{Eq2} and Proposition \ref{Pr25}, we also have
\begin{align}\label{Eq4}
\mE\Lambda\left(|\eta^X_{\tau}(x)-\eta^X_{\sigma}(x)|^{1/2}\right)\preceq_\Lambda\mE\Lambda\left(|\sI^X_{\tau}(x)-\sI^X_{\sigma}(x)|\right).
\end{align}
Now, integrating both sides of \eqref{Eq3} and \eqref{Eq4} with respect to $\mu$ over $\mX$ and by Fubini's theorem, we obtain
$$
\mE[\sI^X_{\tau}-\sI^X_{\sigma}]_\Lambda\preceq_\Lambda\mE [\eta^X_\tau-\eta^X_{\sigma}]_{\Lambda_{1/2}},
$$
and that if $\Lambda\in\sA_2$, then 
$$
\mE [\eta^X_\tau-\eta^X_{\sigma}]_{\Lambda_{1/2}}\preceq_\Lambda\mE[\sI^X_{\tau}-\sI^X_{\sigma}]_\Lambda,
$$
where $\Lambda_{1/2}(t):=\Lambda(t^{1/2})$ still belongs to $\sA_1$. Now, for any $\Phi\in\sA_0$, the desired estimates 
\eqref{AA6} and \eqref{AA61} follow by Corollary \ref{TH1}.
\end{proof}
\br
By Proposition \ref{Pr1}, for $\Lambda\in\sA_1$, \eqref{AA6} can be written as
$$
\mE\left(\Phi\left(\sup_{t\in[0,\tau]}\left\|\sI^X_t\right\|_{\Lambda}\right)\right)\preceq_{\Lambda,\Phi}
\mE\left(\Phi\circ\varphi\circ\phi\left(\Bigg\|\Bigg(\int^\tau_0 \|X_s\|^2_{\ell^2}\dif s\Bigg)^{1/2}\Bigg\|_{\Lambda}\right)\right),
$$
and if $\Lambda\in\sA_2$, then \eqref{AA61}  can be written as
$$
\mE\left(\Phi\circ\phi\circ\varphi\left(\Bigg\|\Bigg(\int^\tau_0 \|X_s\|^2_{\ell^2}\dif s\Bigg)^{1/2}\Bigg\|_{\Lambda}\right)\right)\preceq_{\Lambda,\Phi}
C\mE\left(\Phi\left(\sup_{t\in[0,\tau]}\left\|\sI^X_t\right\|_{\Lambda}\right)\right),
$$
where $\varphi$ and $\phi$ are the same as in \eqref{AA2}.
In particular, if $\Lambda(t)=t^p$ for some $p>1$, then for any $\Phi\in\sA_0$ and stopping time $\tau$,
$$
\mE\left(\Phi\left(\sup_{t\in[0,\tau]}\left\|\sI^X_t\right\|_{p}\right)\right)\asymp_{p,\Phi}
\mE\left(\Phi\left(\Bigg\|\Bigg(\int^\tau_0 \|X_s\|^2_{\ell^2}\dif s\Bigg)^{1/2}\Bigg\|_{p}\right)\right),
$$
where $\|\cdot\|_p$ denotes the usual $L^p$-norm. This type of estimate was established in \cite{NVW, Co-Ve} 
due to the fact that $L^p$-space for $p>1$ is a UMD space.
\er

\end{document}